\documentclass{article}

\usepackage{amsmath}
\usepackage{amsfonts}
\usepackage{amscd}
\usepackage{amssymb}

\input xypic

\newtheorem{defn}{Definition}[section]
\newtheorem{thm}[defn]{Theorem}
\newtheorem{prop}[defn]{Proposition}
\newtheorem{lemma}[defn]{Lemma}

\newtheorem{fact}[defn]{Fact}

\newcommand{\lm}{\ensuremath{\longrightarrow}}

\DeclareMathOperator{\proj}{\mbox{Proj}\,}

\DeclareMathOperator{\Gr}{\mbox{Gr}}

\DeclareMathOperator{\pd}{\mbox{pd}\,}

\DeclareMathOperator{\Quot}{\mbox{Quot}}

\DeclareMathOperator{\s}{\sigma}

\DeclareMathOperator{\Z}{\mathbb{Z}}

\begin{document}

\begin{center}
  \LARGE \textbf{Hilbert Schemes for Quantum Planes are Projective}
\end{center}

\begin{center}
  DANIEL CHAN
\end{center}

\begin{center}
  {\em University of New South Wales}
\end{center}

\begin{center}
e-mail address:{\em danielc@unsw.edu.au}
\end{center}

\begin{abstract}
We show that Hilbert schemes for the quantum plane are projective. We also show that some collections of torsion sheaves are bounded.
\end{abstract}

Throughout this paper, all objects will be defined over a fixed ground field $k$. 

\section{Introduction}

In [AZ2], Artin and Zhang developed the theory of Hilbert schemes in a very general categorical setting. The main motivation for such a generalisation was for applications to the quotient category $\proj R:=\Gr-R/\mbox{tors}$ where $R$ is a non-commutative locally finite connected graded $k$-algebra, $\Gr-R$ is the category of graded $R$-modules and tors is the Serre subcategory of $R_{>0}$-torsion modules. They identified a condition on $R$, namely strong $\chi$, which guarantees that the Hilbert schemes are well-behaved, in particular, they are countable unions of projective schemes. They conjecture [AZ2; conjecture~E5.2] that the Hilbert schemes are moreover projective. 

The main purpose of this note is to verify a special case of this conjecture.

\begin{thm} \label{tmain}
Let $R$ be an AS-regular algebra of dimension 3 generated in degree 1 with 3 quadratic relations. Let $F \in \proj R$ and $h(t)$ be a Hilbert polynomial. Then the Hilbert scheme parametrising quotients of $F$ with Hilbert polynomial $h(t)$ is projective.
\end{thm} 
Recall that $\proj R$ is often referred to as a quantum plane since $R$ is a non-commutative analogue of the commutative polynomial ring in 3 variables. We will also look at moduli schemes of some torsion sheaves on $\proj R$. 

\section{Preliminaries}

For the rest of the paper, we let $R$ denote an AS-regular algebra of dimension 3 that is generated in degree 1 and has 3 quadratic relations. The precise definition of such algebras can be found in [ATV1, section~2]. Instead, we will list in this section, all the important facts about such algebras that we use. We start with
\begin{fact}
The algebra $R$ is a noetherian domain of global dimension three.
\end{fact}

That $R$ is global dimension 3 is part of the definition whilst the fact that $R$ is a noetherian domain is just [ATV1, theorem~8.1] and [ATV2, theorem~3.9].

Let $M$ be a finitely generated graded $R$-module. Since our $R$-modules will always be graded, we will drop the adjective ``graded'' in future. We can associate to $M$ its Hilbert function $h_M(t) = \dim_k M_t$ where $M_t$ is the degree $t$ component. Let $\pi: \Gr R \lm \proj R$ be the quotient functor by the Serre subcategory of $R_{>0}$-torsion modules. Here, by $R_{>0}$-torsion module, we mean modules $M$ where each $m \in M$ is annihilated by a power of the augmentation ideal $R_{>0}$. Recall from [AZ1, section~7] that one can define cohomology of objects in $\proj R$. Now $h_M(t)$ is not a well-defined function of $\pi M$ so we will also use the {\em Hilbert polynomial} of $M$ defined  by 
\[  P_M(t) = \sum_{i\geq 0} (-1)^i \dim H^i(\proj R, \pi M(t))  \] 
where $M(t)$ is the module $M$ with degrees shifted by $t$ so $M(t)_i = M_{i+t}$, and the sum is well-defined by the non-commutative version of Serre's finiteness theorem in [AZ1, theorem~7.4]. Standard cohomology theory shows that $P_M(t) = h_M(t)$ for $t \gg 0$.

\begin{fact} ([ATV1, (2.17)])  
The Hilbert polynomial of $R$ is $P_R(t) =(\begin{smallmatrix} t+2 \\ 2 \end{smallmatrix}) =  \frac{1}{2}(t^2+3t+2)$.
\end{fact}

Since any finitely generated $R$-module $M$ has a finite graded free resolution, $P_M(t)$ is indeed a polynomial. 

\vspace{2mm}

We will also need to introduce the notion of semistability.
We define the {\em slope} of a finitely generated $R$-module $M$ to be 
\[  \mu(M):= \frac{P_M(t)}{e}  \]
where $e$ is the leading coefficient of $P_M(t)$. If $\deg P_M(t)=2$ then $2e$ is just the rank of the module $M$. These polynomial functions may be ordered by their behaviour as $t \lm \infty$ which amounts to the lexicographic order on the coefficients.

\begin{defn}
A finitely generated $R$-module $M$ is said to be {\em semistable} if for every submodule $N \leq M$ we have $\mu(N) \leq \mu(M)$.
\end{defn}
Since $\mu(M)$ depends only on $\pi M$, the definition extends naturally to noetherian objects in $\proj R$. 

As in the commutative case we have,

\begin{lemma}
Consider an exact sequence of $R$-modules $0 \lm M' \lm M \lm M'' \lm 0$ whose Hilbert polynomials all have the same degree. 
\begin{enumerate}
\item $\mu(M)$ is the average (with respect to some positive weights) of $\mu(M'),\mu(M'')$ and so in particular lies in between these two slopes.
\item If $M',M''$ are semistable with the same slope $\mu$ then so is $M$. 
\item For integers $l>0$ and $c$, the module $R(-c)^l$ is semistable. 
\end{enumerate}
\end{lemma}
\textbf{Proof.} We omit the routine proofs of parts i) and ii). Part iii) will follow from ii), once we prove that $R(-c)$ is semistable. Let $N \leq R(-c)$ be non-zero so that $\deg P_N(t) = 2$. Now the ranks of $R(-c)$ and $N$ are both 1 so  $P_N(t) \leq P_{R(-c)}(t)$ implies $\mu(N) \leq \mu(R(-c))$ and $R(-c)$ is indeed semistable.

\section{Proof of main theorem}

In this section, we prove theorem~\ref{tmain}. Artin-Zhang's [AZ2] theory of Hilbert schemes applies to $\proj R$ since $R$ satisfies strong $\chi$ by [AZ2, proposition~C6.10] and the fact that $R$ is strongly noetherian ([ASZ,propositions~4.9(1) and 4.13]) and $R$ has a balanced dualising complex ([Y, theorem~7.18]). We will not need the definition of strong $\chi$ which can be found in [AZ2,C6.8]. 

Let $F$ be a finitely generated $R$-module and $P_1(t)$ be a Hilbert polynomial. Since strong $\chi$ holds, we may apply [AZ2, theorem~E5.1] to see that the Hilbert functor of flat quotients of $\pi F$ with Hilbert polynomial $P_1(t)$ is represented by a Hilbert scheme $\Quot(F,P_1(t))$. Artin and Zhang also give a criterion for when this Hilbert scheme is projective. A diluted version is the following.

\begin{prop} ([AZ2, proposition~E5.10]) \label{pcrit}
Consider the set $\mathcal{N}$ of all $R_{>0}$-torsionfree quotients of $F$ with Hilbert polynomial $P_1(t)$. If the set of Hilbert functions $\{h_N(t)| N \in \mathcal{N} \}$ is also finite, then $\Quot(F,P_1(t))$ is projective.
\end{prop}

To prove theorem~\ref{tmain}, we first reduce to the case $F=R^l$ as follows. Observe that there exists a surjection $\pi R(-c)^l \lm \pi F$ in $\proj R$ for some $c,l \in \Z$. Hence, to show $\Quot(F,P_1(t))$ is projective, it suffices to show $\Quot(R(-c)^l,P_1(t))$ is projective. We may thus assume that $F = R(-c)^l$. Tensoring by $R(c)$ gives an isomorphism between the set of flat families of quotients of $R(-c)^l$ with Hilbert polynomial $P_1(t)$ and the set of flat families of quotients of $R^l$ with Hilbert polynomial $P_1(t+c)$. This induces an isomorphism of Hilbert schemes so we may thus assume henceforth that $F = R^l$. 

Consider a closed point of $\Quot(F,P_1)$ given by an $R_{>0}$-torsionfree module $N$ and a surjection $f: F \lm N$. As in the classical commutative case, we will parametrise such quotients by parametrising the corresponding kernels $M:= \ker f$.  They all have Hilbert polynomial $P(t):= lP_R(t) - P_1(t)$. 

Note that as $N$ is $R_{>0}$-torsionfree, the depth of $N$ is positive. Hence by Auslander-Buchsbaum [Jorg, Theorem~3.2], $\pd N =\mbox{gl.dim}\ R - \mbox{depth}\ N \leq 2$. Consequently, $\pd M \leq 1$ and we can find a minimal resolution of the form 
\[ 0 \lm \bigoplus_{j=1}^m R(-b_j) \lm \bigoplus_{i=1}^n R(-a_i) \lm M \lm 0 .\]
Minimality of the resolution means that the induced maps $R(-b_j) \lm R(-a_i)$ are zero whenever $b_j \leq a_i$. 

The theorem will follow from proposition~\ref{pcrit}, once we show that the possibilities for $a_i,b_j$ are bounded. First note the $a_i\geq 0,\ b_j>0$ since $M \leq R^l$. We may also assume the $a_i, b_j$  are in increasing order. Now the Hilbert polynomial of $M$ is 

\begin{align*}
 P_M(t) & = \sum_{i=1}^n \left( \begin{smallmatrix} t-a_i+2 \\ 2 \end{smallmatrix} \right) - \sum_{j=1}^m \left( \begin{smallmatrix} t-b_j+2 \\ 2 \end{smallmatrix} \right)  \\
  & = \frac{1}{2} (n-m)t^2 + [(\sum b_j - \sum a_i)+ \frac{3}{2}(n-m)]t + [\frac{1}{2}(\sum a_i^2 - \sum b_j^2) + \frac{3}{2}(\sum b_j - \sum a_i) +(n-m)]
\end{align*}
For simplicity, we let $p_i$ denote the coefficient of $t^i$ above. More generally, given two sets of integers $\{a_1,\ldots,a_n\},\{b_1,\ldots,b_m\}$ we consider the following functions on these two sets of integers.
\begin{align*}
p_2(\{a_i\},\{b_j\}) & =  \frac{1}{2}(n-m) \\
p_1(\{a_i\},\{b_j\}) & = (\sum b_j - \sum a_i)+ \frac{3}{2}(n-m) \\
p_0(\{a_i\},\{b_j\}) & = \frac{1}{2}(\sum a_i^2 - \sum b_j^2) + \frac{3}{2}(\sum b_j - \sum a_i) +(n-m)
\end{align*}

If $\pd M = 0$ so that $m=0$, then $\sum a_i = 3p_2 - p_1$ so non-negativity of the $a_i$'s ensures there are only a finite number of possibilities as to what they can be. We assume from now on that $\pd M = 1$. 

For each $r \in \{1,\ldots m\}$ we let $s$ be the largest integer such that $a_s < b_r$. The key numerical constraint we need is 

\begin{lemma}
We have $a_1 + \ldots + a_s \geq b_1 + \ldots + b_r$ and $s>r$.
\end{lemma}
\textbf{Proof.} Let $L$ denote the image of $\oplus^s_{i=1} R(-a_i)$ in $M$ and consider the exact sequence 
\[ 0 \lm K \lm \bigoplus_{i=1}^s R(-a_i) \lm L \lm 0 .\]
Minimality of the free resolution for $M$ shows that $K$ contains $\oplus_{j=1}^r R(-b_j)$ and in fact, we have a direct sum decomposition $K = \oplus^r R(-b_j) \oplus K'$ for some $R$-module $K'$. Now $L$ is a non-zero submodule of $R^l$ so $\deg P_L(t) = 2$ giving the inequality $s>r$. 

We turn our attention to the other inequality. Now $K'$ certainly embeds in  $R^m$ so semistability ensures $\mu(K') \leq t^2+3t+2$. We thus deduce that 
\[ P_{K'}(t) = \frac{1}{2}\rho t^2 + (\frac{3}{2}\rho - \alpha) t + \mbox{const} \]
where $\alpha, \rho \geq 0$. The Hilbert polynomial of $L$ can now be computed as 
\begin{align*}
P_L(t) & = \sum^s \left( \begin{smallmatrix} t-a_i+2 \\ 2 \end{smallmatrix} \right) - \sum^r \left( \begin{smallmatrix} t-b_j+2 \\ 2 \end{smallmatrix} \right)  - P_{K'}(t) \\
& = \frac{1}{2} (s-r)t^2 + [(\sum^r b_j - \sum^s a_i)+ \frac{3}{2}(s-r)]t 
-\frac{1}{2}\rho t^2 + (-\frac{3}{2}\rho + \alpha) t + \mbox{const} \\
& = \frac{1}{2} (s-r-\rho)t^2 + [(\sum^r b_j - \sum^s a_i + \alpha)+ \frac{3}{2}(s-r-\rho)]t + \mbox{const} 
\end{align*}
Semistability gives $\mu(L) \leq t^2+3t+2$ so $\sum^r b_j - \sum^s a_i + \alpha \leq 0$ from which the desired inequality follows. 

\vspace{2mm} 

We continue now with the proof of the theorem. Let $\mathcal{S}$ be the set whose elements consist of the following data
\begin{enumerate}
\item a finite non-decreasing sequence of non-negative integers $0 \leq a_1\leq a_2 \leq \ldots \leq a_n$ for some $n>0$,
\item a finite non-decreasing sequence of positive integers $0 < b_1 \leq b_2 \leq \ldots \leq b_m$ for some $m\geq 0$,
\end{enumerate}
subject to the following condition:

(*) \hspace{1cm} for each $r \in \{1,\ldots,m\}$, if $s(r)$ is the largest integer with $a_{s(r)} < b_r$, then the following inequalities hold:
$$ (*1)\ s(r)>r \quad \mbox{and} \quad (*2)\ a_1 + \ldots + a_{s(r)} \geq b_1 + \ldots + b_r$$
The lemma shows that a minimal resolution of $M$ gives an element of $\mathcal{S}$ with 
$$ p_2(\{a_i\},\{b_j\})=p_2,  p_1(\{a_i\},\{b_j\})=p_1,  p_0(\{a_i\},\{b_j\})=p_0 .$$
It suffices to show that the set of elements of $\mathcal{S}$ with these $p_2,p_1,p_0$ values is finite.

To this end, we perform the following algorithm on the elements $\s=(\{a_i\},\{b_j\})$ of $\mathcal{S}$. Throughout the algorithm, $p_2(\s),p_1(\s)$ will remain unchanged but $p_0(\s)$ strictly increases with each iteration. On each iteration we perform the following steps.

\begin{enumerate}
 \item Pick $s$ maximal such that $a_s < b_1$ (this exists). Reduce both $a_s$ and $b_1$ by 1. 
 \item Relabel indices so that the $a_i$ are in increasing order.
 \item If after relabelling $b_1 = a_s$, then delete $a_s,b_1$ from the sets. Note that if this occurs then some relabelling must have occurred in the previous step and that the deleted values are $a_s = a_{s-1}$ or $a_{s-1}+1$. 
 \item If deletion occurred, relabel the indices so there are no missing $a_i$'s or $b_j$'s, that is, the old $b_{j+1}$ is now $b_j$ for $j\geq 1$, and the old $a_{i+1}$ is now $a_i$ for $i > s$.
\end{enumerate}
The algorithm continues until $\{b_j\}$ is empty.
\begin{lemma}
In the above algorithm, each iteration induces a map $T: \mathcal{S} \lm \mathcal{S}$ satisfying 
$$ p_2(\s) = p_2(T\s), p_1(\s) = p_1(T\s), p_0(\s) < p_0(T\s) .$$
\end{lemma}
\textbf{Proof.} One verifies easily that $p_2(\s),p_1(\s)$ stay fixed. In step~i) however, $p_0(\s)$ strictly increases whilst staying fixed throughout all the other steps. Note $a_1 + \ldots + a_s \geq b_1 >0$ so $a_s >0$. Since only $a_s$ decreases in each iteration, we see the $a_i$'s remain non-negative. This shows that if $b_1 = 1$, it will be deleted in step iii), so the $b_j$'s will always remain positive. Finally, we check preservaton of the inequalities in (*). If no deletion occurred in step iii), then the function $s(r)$ remains unchanged and both sides of the inequality (*2) are reduced by 1. If deletion did occur, then $s(r)$ changes to $s(r+1)-1$ so (*1) still holds and both sides of (*2) are reduced by $a_s=b_1$. This completes the proof of the lemma. \vspace{2mm}

We now complete the proof of the theorem by running the algorithm on the element of $\mathcal{S}$ given by the minimal resolution of $M$. We note that the algorithm terminates since the $b_j$'s are positive. When the algorithm terminates we end up with an element $\s \in \mathcal{S}$ given by integers 
\[ 0 \leq a_1 \leq a_2 \leq \ldots \leq a_{n-m}  .\]
Furthermore, their sum is $3p_2-p_1$ so there are only a finite number of possibilities for these final values of $a_i$ and hence, only a finite number of possibilities for the final value of $p_0(\s)$. Comparing this with the original value of $p_0$ we see that at most $p_0(\s) - p_0$ iterations could have occurred in our algorithm. We now run the algorithm in reverse to examine the possible original values for $a_i,b_j$. On each iteration, if values of $a_s,b_1$ were deleted then by the remark in step iii), they are bounded by the other values. The only other possible change to the values is in step i) so there is an explicit upper bound on the original values of $a_i,b_j$. This completes the proof of the theorem.

\vspace{5mm}

\section{Moduli of torsion sheaves}

In the commutative theory, Hilbert schemes are useful in the construction of many moduli schemes. We show how our methods can similarly be applied to study the moduli of some torsion sheaves on $\proj R$. 

Recall that given a collection of (isomorphism classes of) coherent sheaves $\mathcal{M}$ on a scheme $X$, we say that $\mathcal{M}$ is {\em bounded} if there exists an algebraic variety $S$ and a flat family of sheaves $\mathcal{N}$ on $X$ parametrised by $S$, such that every sheaf in $\mathcal{M}$ occurs as some closed fibre $\mathcal{N}\otimes_S k(s), s \in S$. This definition extends naturally to $X = \proj R$.

The sheaves on $\proj R$ we wish to study in this section can be described as follows. Let $M \in \Gr R$ be a finitely generated module with Hilbert polynomial $P(t) = et+d$ for some $e,d \in \Z$. We say $M$ is {\em pure dimension one} if furthermore $\pd M = 1$. The corresponding object $\pi M \in \proj R$ is called a {\em pure 1-dimensional coherent sheaf on $\proj R$}. 

\begin{thm}
The collection of all pure 1 dimensional semistable sheaves on $\proj R$ of Hilbert polynomial $et+d$ is bounded.
\end{thm}
\textbf{Proof.} 
Consider the minimal resolution 
\[ 0 \lm \bigoplus_{j=1}^n R(-b_j) \lm \bigoplus_{i=1}^n R(-a_i) \lm M \lm 0 .\]
We assume that both the $a_i$'s and $b_j$'s are in increasing order. Since Hilbert schemes are projective, it suffices to show that the possibilities for the $a_i$'s are finite.


We need to prove the key
\begin{lemma}
For $m \in \{1,\ldots n-1\}$ we have $b_m > a_{m+1}$.
\end{lemma}
\textbf{Proof.} First observe that $b_m > a_m$ for otherwise, by minimality of the resolution we find that $\oplus_{j=1}^m R(-b_j)$ must embed in $\oplus_{i=1}^{m-1} R(-a_i)$. Rank considerations show this is impossible so certainly $b_m > a_m$. Suppose now that $b_m \leq a_{m+1}$. This time minimality of the resolution shows that $\oplus^m R(-b_j)$ embeds into $\oplus^m R(-a_i)$. Furthermore, if $M'$ denotes the image of $\oplus^m R(-a_i)$ in $M$ then rank considerations show that we have an exact sequence 
\[ 0 \lm \bigoplus_{j=1}^m R(-b_j) \lm \bigoplus_{i=1}^m R(-a_i) \lm M' \lm 0 .\]
We seek to show that $\mu(M') > \mu(M)$ and so obtain a contradiction. Now the Hilbert polynomial of $M$ is 
\[ P_M(t) = \sum_{i=1}^n(b_i-a_i)t + \frac{1}{2}\sum_{i=1}^n(b_i-a_i)(3-a_i-b_i)  \]
and a similar formula holds for $P_{M'}(t)$. Hence if $\mu(M)=t+c$ then $c$ is given by the weighted average of $\frac{1}{2}(3-a_i-b_i)$ where the weights are $b_i-a_i>0$. Now the values of $3-a_i-b_i$ are non-increasing and we have a strict inequality 
\[ 3-a_m-b_m > 3-a_{m+1}-b_{m+1}    \]
by our hypothesis. The inequality $\mu(M') > \mu(M)$ and hence lemma now follow from the following elementary
\begin{fact}
Let $w_1,\ldots,w_n>0$ be positive weights and $x_1\geq x_2 \geq \ldots \geq x_n$ be a non-increasing sequence of numbers. Then for any $m \in \{1,\ldots,n\}$ we have  
$$\frac{w_1x_1 + \ldots + w_m x_m}{w_1 + \ldots + w_m} \geq \frac{w_1x_1 + \ldots + w_n x_n}{w_1 + \ldots + w_n}$$
with equality iff either $m=n$ or $x_1=x_n$. 
\end{fact} 

\vspace{2mm}

We continue the proof of the theorem. As already noted in the proof of the lemma, the $b_i - a_i$ are positive. Also, they sum to $e$,  so the lemma shows that all the values of $a_i,b_j$ lie in the interval $[a_1,a_1+ e]$. Since $d = \frac{1}{2}\sum(b_i-a_i)(3-a_i-b_i)$ is fixed, $|a_1|$ cannot be arbitrarily large. This completes the proof of the theorem. 

\vspace{1cm}
\textbf{\large{References}}

\begin{itemize}
  \item[{[AZ1]}]  M. Artin, J. Zhang, ``Noncommutative projective schemes'', {\em Adv Math.}, \textbf{109} (1994) p.228-87
  \item[{[AZ2]}]  M. Artin, J. Zhang, ``Abstract Hilbert Schemes'' {\em J. Algebras  and Repr. Theory}, \textbf{4} (2001) p.305-94
  \item[{[ASZ]}] M. Artin, L. Small,  J. Zhang, ``Generic flatness for strongly noetherian rings'' {\em J. Algebra} \textbf{221} (1999) p.579-610
  \item[{[ATV1]}]  M. Artin, J. Tate, M. van den Bergh, ``Some algebras associated to automorphisms of elliptic curves'', {\em The Grothendieck Festschrift, Vol. 1} p.33-85, Progress in Math., Birkhauser, Boston, (1990)
  \item[{[ATV]}]  M. Artin, J. Tate, M. van den Bergh, ``Modules over regular algebras of dimension 3'', {\em Invent. Math.} \textbf{106} (1991) p.335-88
  \item[{[Jorg]}] P. Jorgensen, ``Non-commutative graded homological identities'', {\em J. London Math. Soc.} \textbf{57} (1998) p.336-50
  \item[{[Y]}] A. Yekutieli, ``Dualizing complexes over noncommutative graded algebras'', {\em J.. Algebra} \textbf{153} (1992) p.41-84
\end{itemize}

\end{document}